\documentclass[11pt,a4paper]{article}
\usepackage{natbib}
\usepackage[resetlabels]{multibib}
\usepackage[utf8]{inputenc}
\newcites{latex}{References}
\usepackage{amsmath}
\usepackage{leftidx}
\usepackage{amsthm}
\usepackage{todonotes}
\usepackage{latexsym}
\usepackage{amssymb}
\usepackage{bbm}
\usepackage[english]{babel}
\usepackage{array}
\usepackage{rotating}
\usepackage{graphicx, color}
\usepackage{multirow}
\usepackage{pslatex}
\usepackage{lscape}
\usepackage[flushmargin, hang]{footmisc}
\usepackage{fancyhdr}
\usepackage{dsfont}
\usepackage{float}
\usepackage{arydshln}
\usepackage{enumitem}
\usepackage{setspace}
\usepackage[T1]{fontenc}
\usepackage{ae, aecompl}
\usepackage{sectsty}
\usepackage{tikz}
\usepackage{comment}
\usepackage[labelfont= bf,sf,small]{caption}
\usepackage[breaklinks, pdfstartview= FitH, colorlinks=true,
linkcolor=blue, citecolor=blue, urlcolor=black]{hyperref}
\usepackage{marginnote}
\usepackage{nameref}
\usepackage{adjustbox}

\allsectionsfont{\sffamily}
\newtheorem{Theorem}{Theorem}[section]
\newtheorem{Lemma}[Theorem]{Lemma}

\newtheorem{Definition}[Theorem]{Definition}

\theoremstyle{nonumberplain}
\newtheorem{Proof}{Proof}

\newcommand{\beq}{\begin{equation}}
\newcommand{\eeq}{\end{equation}}
\newcommand{\beqs}{\begin{equation*}}
\newcommand{\eeqs}{\end{equation*}}
\newcommand{\beqn}{\begin{eqnarray}}
\newcommand{\eeqn}{\end{eqnarray}}
\newcommand{\bs}{\begin{eqnarray*}}
\newcommand{\es}{\end{eqnarray*}}
\newcommand{\bi}{\begin{itemize}}
\newcommand{\ei}{\end{itemize}}
\newcommand{\be}{\begin{enumerate}}
\newcommand{\ee}{\end{enumerate}}
\newcommand{\bd}{\begin{description}}
\newcommand{\ed}{\end{description}}

\long\def\symbolfootnote[#1]#2{\begingroup%
\def\thefootnote{\fnsymbol{footnote}}\footnote[#1]{#2}\endgroup}

\title{Roth's Theorem implies a Weakened Version of the ABC Conjecture for Special Cases\thanks{The authors are grateful to Enrico Bombieri, Michel Waldschmidt, Maxie D. Schmidt and Joshua Lampert for helpful comments.}}

\author{Philipp Sibbertsen, Karsten Müller, Timm Lampert, Michael Taktikos}
\begin{document}

\maketitle

\begin{abstract}
Enrico Bombieri proved that the ABC Conjecture implies Roth's theorem in 1994. This paper concerns the other direction. In making use of Bombieri's and Van der Poorten's explicit formula for the coefficients of the regular continued fractions of algebraic numbers, we prove that Roth's theorem implies a weakened non-effective version of the ABC Conjecture in certain cases relating to roots.
\end{abstract}


\section{Introduction}

The ABC Conjecture is one of the most famous open problems in number theory. Masser and Oesterle conjectured it in \cite{Masser} and \cite {Oesterle}. We investigate connections to Roth's theorem, which was proven by \cite{Roth}. \cite{Bombieri} proved that the ABC Conjecture implies Roth's theorem. We consider the other direction under weakened assumptions for special cases.

\section{The ABC Conjecture}

To formulate the ABC Conjecture we need the following definition:

\begin{Definition}For a positive integer $a$, $rad(a)$ is the product of the distinct prime factors of $a$. 
\end{Definition}

\begin{Definition}[ABC Conjecture]
For every positive real number $\varepsilon$, there exists a constant $K_{\varepsilon}$
such that for all triples $(a,b,c)$ of coprime positive integers,
with $a + b = c$:
\[c < K_{\varepsilon} \cdot rad(abc)^{1+\varepsilon}.\]
\end{Definition}

Note that this version of the ABC Conjecture is not effective. This is so because only the existence of a constant $K_{\varepsilon}$ is claimed without specifying how to calculate $K_{\varepsilon}$ for given $\epsilon$. 

Explicit effective versions which specify $K_{\varepsilon}$, however, are used in proofs
related to equations of form $a^{m} + b^{n} = c^{k}$. In such proofs, the use of an effective ABC Conjecture can ignore prime factorisation and, thus, $rad(abc)$. Instead, it suffices to refer to $abc$ without exponents. This motivates our definition of a ``weakened effective version of ABC'' (WEABC).

\begin{Definition}[WEABC]
 For any fixed $\epsilon > 0$ there is a constant $K_{\varepsilon}$ such that for all triples $(a,b,c)$ of coprime  positive integers and $m,n,k \in \mathbb{N}$, with $a^{m} + b^{n} = c^{k}$:
\[c^{k} < K_{\varepsilon} \cdot (abc)^{1+\varepsilon}.\]
\end{Definition}

This definition differs from the literature where effective versions of the ABC Conjecture for a fixed $\varepsilon > 0$ are called ``weak ABC Conjecture'', e.g., by \citep[p. 403, example 12.2.3]{Gubler}.

Fermat's last theorem follows immediately from WEABC for $n\geq 6$ \citep[p.403, example 12.2.4]{Gubler}:

\begin{Theorem}The WEABC with $\varepsilon = 1$ and $K_{\varepsilon} =1$ implies Fermat's Theorem for $n\geq 6$.
\end{Theorem}

\begin{Proof}
Let $x,y$ and $z$ be coprime positive integers satisfying $x^n + y^n = z^n$.
With $a = x^n$, $b = y^n$, $c = z^n$ we get $z^n \leq rad(xyz)^2 \leq (xyz)^2 \leq z^6$.
As $z > 1$ the desired inequality $n \leq 5$ follows.\qed
\end{Proof}

Likewise, Catalan's conjecture for large $n$ and $m$ can be proven from WEABC.

Furthermore, theorems that have been proven independently of the ABC Conjecture and that could be proven by the ABC Conjecture like \cite{Wiles} proof of Fermat's conjecture, WEABC suffices to prove those theorems. 
So it is enough to get rid off the exponent $n$.

If we do not refer to a specified $K_{\varepsilon}$, we speak of a ``weakened non-effective ABC Conjecture'' (WNEABC). 

\begin{Definition}[WNEABC]
 For every $\epsilon > 0$, there exists a constant $K_{\varepsilon}$ for all triples $(a,b,c)$ of coprime  positive integers and $m,n,k \in \mathbb{N}$, with $a^{m} + b^{n} = c^{k}$:
\[c^{k} < K_{\varepsilon} \cdot (abc)^{1+\varepsilon}.\]
\end{Definition}

Such a conjecture can be proven from Roth's theorem because it is not necessary to find estimates for the prime numbers, which is extremely difficult. A generalised version of Ridout's theorem would even prove the full ABC Conjecture for these cases, but this seems to be very hard.

\begin{Theorem}Roth's Theorem:
Let $a$ be an algebraic number. Then for every $\varepsilon > 0$, there exists a constant $K_{\varepsilon}$ such that for all positive integers $p$ and $q$:
\begin{equation*}|a - \frac{p}{q}| > \frac{K_{\varepsilon}}{q^{2+\varepsilon}}.\end{equation*}
\end{Theorem}
\begin{Proof}
See \cite{Roth}.\qed
\end{Proof}

Roth's theorem is not effective. Therefore, only non-effective ABC cases can follow from it.
Roth's theorem can probably be sharpened. It is even an open problem whether the coefficients of the regular continued fraction of algebraic numbers of degree $> 2$ are bounded or not.
If Roth's theorem can be sharpened then the  WNEABC  can be sharpened as well.

Counterexamples for the  ABC Conjecture with $\varepsilon=0$ are known. But these counterexamples do not result from algebraic equations. 
Therefore, the coefficients of the regular continued fractions of algebraic numbers of degree greater than 2 can nonetheless be bounded.

\begin{Definition}A triple of coprime integers $a+b=c$ is an ABC hit, if $rad(abc) < c$ holds and the quality of a triple is given by the first $quality(a,b,c) = \frac{ln(c)}{ln(rad(abc))}$. As $c>1$ it is clear that abc is larger than 1.
\end{Definition}

$2 + 23^5 = 9^5 \cdot 109$  is the ABC hit with the highest quality so far found. Its quality is approximately 1.62991.
This hit results from the 3rd convergent to $\sqrt [5]{109}$ with regular continued fraction $[2;1,1,4,77733,\ldots]$. The fraction $\frac{23}{9}$ is a good approximation because 77733 is very large and very near the absolute Liouville bound.

Thus, ABC hits are connected with diophantine approximations. We want to go further in this direction. Each number of the form $\sqrt [s]{k}$ leads to a set of $abc$ equations due to the convergents of the root.
However, the resulting $a,b$ and $c$ are not necessarily coprime. This is not a restriction because the ABC Conjecture implies Roth's theorem, as proven by \cite{Bombieri}. We will show by Theorem \ref{maintheorem} that WNEABC is also valid for these equations even when the resulting $a,b$ and $c$ are not coprime.

Each root $\sqrt [s]{k}$ is directly connected with a set of $abc$ equations.
\begin{Definition}Resulting equations for $\sqrt [s]{k}$:=\label{resultingequ}
Let $\frac{p_n}{q_n}$ be a convergent of  $\sqrt [s]{k}$ and $d_n = k q_n^s - p_n^s$.
If $d_n > 0$ then the resulting equation is
$d_n + p_n^s = k q_n^s$, otherwise it is
$-d_n + k q_n^s = p_n^s$.
\end{Definition}

$p_n$ and $q_n$ are always coprime, but $k$ and $d_n$ may have common prime factors with the other factors of the equations.

From the ABC record hits known until now most hits have resulting equations, cf. Table \ref{tablehits} for the most famous ones. It is remarkable that several of the best ABC hits have very small $a,b,c$.

\begin{table}
\begin{tabular}{|c|c|c|}\hline
{\bf abc Equation} & {\bf Root} & {\bf Founder}\\ \hline
 2 + $23^5$ = $9^5\cdot 109$ &  $\sqrt[5]{109}$   & Reyssat \\ \hline
  $5^{56} \cdot 245983 = 2^{15} 3^{77} \cdot 11 \cdot 173 + 2543^4 \cdot 182587 \cdot 2802983$ &  $\sqrt[56]{245983}$  &  Bonse\\ \hline
$5^4 \cdot 7 = 2 3^7 + 1$ &  $\sqrt[4]{7}$ & De Weger \\\hline
\end{tabular} 
\caption{Particular cases of abc Equations and \label{tablehits} Roots}
\end{table}

For our WEABC, a 2nd quality results in these cases, which only considers the approximation speed and not the prime factorisation of the convergents:

\begin{Definition}
The 2nd quality(a,b,c) equals $\frac{ln({q_n}^{s} k)}{ln(d_n q_n k p_n)}$ if $d_n>0$, and
$\frac{ln({p_n}^{s})}{ln(-d_n q_n k p_n)}$ if $d_n<0$.
\end{Definition}

The 2nd quality is always smaller or equal to the 1st quality.
There are already ABC hits regarding the 2nd quality alone, in the sense that the 2nd quality is larger than 1. The WNEABC is known to be equivalent to the conjecture
that the 2nd quality has Limes Superior $\leq 1$ for the resulting equations of the convergents of a root.

The main aim of this paper is to show how to prove the following theorem:

\begin{Theorem}\label{maintheorem}
The  WNEABC is valid for all resulting equations for $\sqrt [s]{k}$.
\end{Theorem}

\section{Proof of Theorem \ref{maintheorem}}

$d_n, p_n$ and $q_n$ need not to be coprime. But often they will be coprime. So not all cases are ``real ABC cases'', but the inequality of the WNEABC is always valid.

We first illustrate how to prove Theorem \ref{maintheorem} by proving it for quadratic irrationals $\sqrt [2]{k}$. In these cases all $d_n$ are 1 or -1. Furthermore, the coefficients $b_n$ are periodic and therefore bounded. But this latter fact is not even needed.

\begin{Proof}
Case 1: $d_{n} = 1$. According Definition \ref{resultingequ}, $1 + {p_n}^2 = k\cdot {q_n}^2$.
So for the 2nd quality the following results:

\begin{center}
\setlength{\arraycolsep}{0.75mm}
$\begin{array}{lllllllll}
\mbox{2nd quality} &
= & \frac{ln({q_n}^2 \cdot k)}{ln(p_n \cdot k \cdot q_n)} &
= & \frac{2 \cdot ln(q_n) + ln(k)}{ln(p_n) + ln(k) + ln(q_n)}
& = & \frac{2 + \frac {ln(k)} {ln(q_n)}}{1 + \frac {ln(k)}{ln(q_n)} + \frac{ln(p_n)}{ln(q_n}}
& = & \frac{2 + \frac {ln(k)} {ln(q_n)}}{1 + \frac {ln(k)}{ln(q_n)} + \frac{ln(\frac{p_n}{q_n}) + \ln(q_n)} {ln(q_n)}}
\end{array}$
\end{center}

and this converges to 1 as all terms are bounded apart from the unbounded $ln(q_n)$.

Case 2: $d_n = -1$ results similarly in
\begin{center}
$\begin{array}{c}
\frac{2} {2 + \frac {ln(k)}{ln(p_n)} + \frac{ln(\frac{q_n}{p_n})} {ln(p_n)}}
\end{array}$
\end{center}

which converges to 1 as well. \qed
\end{Proof}

When there are laws for the prime factors of the denominator of the quality they can be taken into account for specifying a 3rd quality. For example, for $\sqrt [2]{2}$ with $d_n<0$ this can be done with the prime factor 2 since $q_n$ is always even.

This can be seen from the recurrence relations of the coefficients of the regular continued fraction.
$q_1 = 2$ is even and $q_{n+2} = 2 q_{n+1} + q_n = 5 q_n + 2 q_{n-1}$, therefore
all $q_n$ are even, if $n$ is odd.
This leads to the following specification of the 3rd quality of this case:

\begin{center}
$\begin{array}{lllllll}
\mbox{3rd quality} & = & \frac{ln({p_n}^2)}{ln(p_n q_n)} &
=  & \frac{2} {2 + \frac{ln(\frac{q_n}{p_n})} {ln(p_n)}} &
= &\frac{2} {2 + \frac{ln(q_n)}{ln(p_n)}) - 1}.
\end{array}$
\end{center}

The 3rd quality converges monotonously decreasing to 1 because the denominator is monotonously increasing and so the maximum of the 3rd quality is approximately 1.226 at the resulting equation 1 + 8 = 9, which is a real ABC equation. In this case all resulting equations are real ABC equations as $q_n$ and $p_n$ are coprime and $q_n$ is even so that $p_n$ must be odd and the prime factor $k=2$ cannot interfere.

The 3rd quality is always larger than the 2nd quality. 

 The laws for the prime factorisation of the convergents, however, are not sufficient to get the full ABC inequality with $rad(abc)$. For example, it is known that from the resulting Pell numbers in the denominators of the convergents of $\sqrt [2]{2}$
 the only squares, cubes, or any higher power of an integer are 0, 1, and $169 = 13^2$ (for more details see \cite{Cohn}). An analogous result for Fibonacci numbers has been obtained by \cite{Bugeaud}.

The main tool for the following proofs is the explicit formula for the coefficients of the regular continued fractions of algebraic numbers by \citep[p.151, Theorem 3, formula (13)]{Poorten}, in which${f(x)}$ is the minimal polynomial of the algebraic number. The formula connects $d_n$ and $b_n$ so that WNEABC can be proven from Roth's theorem, which has no direct connection to $d_n$. In the main term of the formula, 
$\frac{f'(x)} {f(x)}$ leads to $d_n$ in the denominator. This holds because $\frac{f'(x)} {f(x)}$  results in $\frac{f'(\frac{p_n}{q_n})} {f(\frac{p_n}{q_n})} = \frac {q_n^s s p_n^{s-1}} {d_n q_n^{s-1}}$ for the convergents  with $f(x) = x^s - k$ being the minimal polynomial of the $s$th root of k.

For the 3rd degree the error term of the formula \citep[p.151, Theorem 3, formula (13)]{Poorten} can be ignored due to the effective bound by Liouville's theorem. For higher order it can always be controlled, when $n$ is large enough, which suffices here to prove WNEABC.

To illustrate the use of Bombieri's and van der Poorten's formula, we first provide a criterion for an ABC hit regarding the 2nd quality of $\sqrt [3]{2}$.

\begin{Theorem}
If $b_{n+1} > 6$ then the resulting equation from $\frac{p_n}{q_n}$ is an ABC hit even regarding the 2nd quality.
\end{Theorem}

\begin{Proof}
Let $d_n > 0$. Then the ABC hit condition for the 2nd quality is

\begin{center}
$\frac{\ln(p_n^3)} {\ln(2 \cdot p_n \cdot q_n \cdot d_n)} > 1.$
\end{center}

Since the logarithm is monotonously increasing, the last equation is equivalent to

\begin{center}
$p_n^3 > 2 \cdot p_n q_n d_n.$
\end{center}

The previous equation, in turn, is equivalent to 

\begin{center}
$\frac{p_n^2}{2 q_n} > d_n.$
\end{center}

From \citep[p.151, Theorem 3, formula (13)]{Poorten}, we have that

\begin{center}
$b_{n+1} < \frac{3 p_n^2}{q_n d_n}.$
\end{center}

Therefore, we find that

\begin{center}
$d_n < \frac{3 p_n^2}{q_n b_{n+1}}.$
\end{center}

Hence, the ABC 2nd quality hit condition $p_n^3 > 2 p_n q_n d_n$ is satisfied for $b_{n+1} > 6$
because the inequality $\frac{3 \cdot p_n^2}{q_n b_{n+1}} < \frac{p_n^2}{2 \cdot q_n}$ holds in this case.\qed
\end{Proof}

We now illustrate the proof of Theorem  \ref{maintheorem} for $\sqrt [3]{2}$. We start with the case that the coefficients $b_n$ are bounded. In this case, $\varepsilon=0$ can even be chosen in WNEABC. 

\begin{Theorem}
The  WNEABC is valid for all resulting equations for $\sqrt[3]{2}$.
\end{Theorem}

\begin{Proof}\label{specialproof}
The explicit Bombieri / van der Poorten formula for $\sqrt [3]{2}$ for the coefficients of the regular continued fraction is

\begin{center}
$b_{n+1} = \lfloor (\frac{{(-1)}^{n+1}}{q_n} \frac{3 p_n^2}{-d_n} - \frac{q_{n-1}}{q_n})\rfloor.$
\end{center}

Let the natural number $K$ be the bound of $b_n$. By the formula of \citep[p.151, Theorem 3, formula (13)]{Poorten} the following follows from $b_{n+1} < K$:

\begin{center}
$\begin{array}{ccc}
\frac{3 {p_n}^2}{q_n |d_n|} < (K+2) & \Leftrightarrow & \frac{3 {p_n}^2} {q_n (K+2)} < |d_n|.  
\end{array}$
\end{center}

To get WNEABC with $\varepsilon=0$ a new constant $L$ must be found such that
\begin{center}
$2 q_n^3 < 2 \cdot L \cdot p_n q_n |d_n|$ if $d_n < 0$;\\
$p_n^3 < 2 \cdot L \cdot  p_n q_n d_n$ if $d_n > 0$.
\end{center}
This leads to 
\begin{center}
$\frac{q_n^2} {p_n L} < |d_n|$
\end{center}
and
\begin{center}
$\frac{p_n^2} {2 q_n L} < d_n.$
\end{center}

Thus $L$ must be chosen to satisfy both
\begin{center}
$\frac{q_n^2} {p_n L} < \frac{3 {p_n}^2} {q_n (K+2)}$
\end{center}
and
\begin{center}
$\frac{p_n^2} {2 q_n L} < \frac{3 {p_n}^2} {q_n (K+2)}.$
\end{center}

This is possible because $\frac{p_n}{q_n}$ converges to $\sqrt [3]{2}$.

It follows that WNEABC follows when $\varepsilon=0$.

To prove  that WNABC hold when $\varepsilon>0$, we need the following Lemma:

\begin{Lemma}Roth's Theorem implies that $\frac{p_n^{\varepsilon}}{b_{n+1}}$ is bounded away from zero for every $\varepsilon > 0$.
\end{Lemma}

\begin{Proof}
This Lemma is proven by the following equations:

\setlength{\arraycolsep}{0.75mm}
\begin{center}
$\begin{array}{llllllll}
0 < K < & |a - \frac{p_n}{q_n}| {q_n}^{2+\varepsilon} & < & |\frac{p_{n+1}}{q_{n+1}} - \frac{p_n}{q_n}| {q_n}^{2+\varepsilon} & = & \frac {1}{q_n q_{n+1}} {q_n}^{2+\varepsilon}\\
= & \frac {q_n} {q_{n+1}} {q_n}^{\varepsilon}
& = & \frac {q_n} {b_{n+1} q_n+q_{n-1}} {q_n}^{\varepsilon} & =  & \frac {1} {b_{n+1} + \frac{q_{n-1}}{q_n}} {q_n}^{\varepsilon} & \leq &\frac {{q_n}^{\varepsilon}} {b_{n+1}}
\end{array}$
\end{center}\qed
\end{Proof}

The weakened non-effective ABC condition is
\begin{center}
$\begin{array}{lllll}
K_{\varepsilon} (2 \cdot d_n p_n q_n)^{1+\varepsilon} & \geq &
K_{\varepsilon} ({\frac{6 p_n^3}{b_{n+1}+2}})^{1+\varepsilon} & \geq &
K_{\varepsilon} \cdot 6^{1+\varepsilon} p_n^3 (\frac{p_n^\varepsilon}{b_{n+1} + 2}) (\frac{p_n^2}{b_{n+1} + 2})^\varepsilon
\end{array}$
\end{center}

The other case $d_n > 0$ can be proven analogously.

Because $(\frac{p_n^\varepsilon}{b_{n+1} + 2})$ is bounded away from zero by Roth's theorem and $(\frac{p_n^2}{b_{n+1} + 2})$ is bounded away from zero by the formula in \citep[p.151, Theorem 3, formula (13)]{Poorten}, $b_{n+1}$ behaves like $\frac{p_n}{d_n}$. 
Therefore, $K_\varepsilon > 0$ can be chosen, so that the  weakened non-effective ABC condition can be satisfied, which concludes the proof.\qed
\end{Proof}

To apply this proof to other roots of the form $\sqrt[s]{k}$, all that has to be done is to exchange numbers. Therefore, Proof \ref{specialproof} for  $\sqrt[3]{2}$ generalises to all roots $\sqrt[s]{k}$, which proves Theorem \ref{maintheorem}. However, it has to be noticed that each root has its own
 $K_\varepsilon$, which is ensured by Roth's theorem.
Note that for higher degree in \citep[p.151, Theorem 3, formula (13)]{Poorten} formula (13) is valid and the error term can be ignored for $n$ large enough. So, our case can always be solved.


{}


\begin{thebibliography}{}

\bibitem[Bombieri(1994)]{Bombieri} Bombieri, E.(1994): Roth’s Theorem and the abc-Conjecture, preprint ETH Zürich.

\bibitem [Bombieri $\&$ Gubler(2006)]{Gubler} Bombieri, E. and Gubler, W. (2006): Heights in Diophantine geometry, New Mathematical Monographs 4. Cambridge: Cambridge University Press.

\bibitem[Bombieri $\&$ van der Poorten(1975)]{Poorten} Bombieri, E. and van der Poorten, A. (1975): ``Continued Fractions of Algebraic Numbers'', in: Baker (ed.), \emph{Transcendental Number Theory},
Cambridge University Press, Cambridge, 137-155.

\bibitem[Bugeaud(2006)]{Bugeaud} Bugeaud, Y; Mignotte, M; Siksek, S (2006): "Classical and modular approaches to exponential Diophantine equations. I. Fibonacci and Lucas perfect powers", Ann. Math., 2 (163): 969–1018

\bibitem[Cohn(1996)]{Cohn} Cohn, J. H. E. (1996): \emph{Perfect Pell powers}, Glasgow Mathematical Journal. 38 (1): 19–20.

\bibitem[Masser(1985)]{Masser} Masser, D. W. (1985): "Open problems". In Chen, W. W. L. (ed.). Proceedings of the Symposium on Analytic Number Theory. London: Imperial College.

\bibitem[Mihailescu(2004)]{Mihailescu}Mihăilescu, P. (2004): "Primary Cyclotomic Units and a Proof of Catalan's Conjecture", J. Reine Angew. Math., 2004 (572): 167–195, 

\bibitem[Oesterle(1988)]{Oesterle} Oesterlé, Joseph (1988), "Nouvelles approches du "théorème" de Fermat", Astérisque, Séminaire Bourbaki exp 694 (161): 165–186.

\bibitem[Roth(1955)]{Roth} Roth, K. (1955): ``Rational Approximations to Algebraic Numbers and Corrigendum'', \emph{Mathematika} 2, 1-20 and 168.

\bibitem [Frankenhuysen(1999)]{Frankenhuysen} Van Frankenhuysen, M. (1999): \emph{The ABC Conjecture implies Roth's theorem and Mordell's conjecture}, Mat. Contemp. 16 (1999), 45-72.

\bibitem[Frankenhuysen(2000)]{Frankenhuysen1} Van Frankenhuysen, M.: \emph{A lower bound in the ABC Conjecture}, J. Number Theory 82, p. 91–95

\bibitem[Frankenhuysen(1995)]{Frankenhuysen2} Van Frankenhuysen, M.: \emph{Hyperbolic spaces and the ABC Conjecture} Dissertation Nijmegen 1995

\bibitem[Waldschmidt(2008)]{Waldschmidt} Waldschmidt, M. (2008): \emph{Introduction to Diophantine Methods}

\bibitem[Wiles(1995)]{Wiles} Wiles, A.J. (1995): \emph{Modular elliptic curves and Fermat’s last theorem}, Ann. of Math. (2) 141, no. 3, 443–551.
\end{thebibliography}
\end{document}